%
%%%%%%%%%%%%%%%%%%%%%%%%%%%%%%%%
%
% These are the macroes
%
%%%%%%%%%%%%%%%%%%%%%%%%%%%%%%%%
\newif\ifsect\newif\iffinal
\secttrue\finalfalse
\def\thm #1: #2{\medbreak\noindent{\bf #1:}\if(#2\thmp\else\thmn#2\fi}
\def\thmp #1) { (#1)\thmn{}}
\def\thmn#1#2\par{\enspace{\sl #1#2}\par
        \ifdim\lastskip<\medskipamount \removelastskip\penalty 55\medskip\fi}
\def\square{{\msam\char"03}}
\def\qedn{\thinspace\null\nobreak\hfill\square\par\medbreak}
\def\pf{\ifdim\lastskip<\smallskipamount \removelastskip\smallskip\fi
        \noindent{\sl Proof\/}:\enspace}
\def\itm#1{\item{\rm #1}\ignorespaces}

\def\bar#1{\overline{#1}}
%
%\input pdfcolor
%
%\input boxedeps
%\input boxedeps.cfg
%\HideDisplacementBoxes
%%
%\def\Figuraeps #1 (#2){\message{Figura #1}
%	\midinsert      
%	\centerline{\BoxedEPSF{#2.eps}}
%	\bigskip
%	\centerline{\bf Figure~#1}
%	\endinsert}
%%
%\def\Figurascaledeps #1 (#2)(#3){\message{Figura #1}
%	\midinsert      
%	\centerline{\BoxedEPSF{#2.eps scaled #3}}
%	\bigskip
%	\centerline{\bf Figure~#1}
%	\endinsert}
%
%
\newcount\parano
\newcount\eqnumbo
\newcount\thmno
\newcount\versiono
%\newcount\remno
\newbox\notaautore
\def\neweqt#1$${\xdef #1{(\number\parano.\number\eqnumbo)}
    \eqno #1$$
    \global \advance \eqnumbo by 1}
\def\newrem #1\par{\global \advance \thmno by 1
    \medbreak
{\bf Remark \the\parano.\the\thmno:}\enspace #1\par
\ifdim\lastskip<\medskipamount \removelastskip\penalty 55\medskip\fi}
\def\newex #1\par{\global \advance \thmno by 1
    \medbreak
{\it Example \the\parano.\the\thmno:}\enspace #1\par
\ifdim\lastskip<\medskipamount \removelastskip\penalty 55\medskip\fi}
\def\newthmt#1 #2: #3{ \global \advance \thmno by 1\xdef #2{\number\parano.\number\thmno}
    \medbreak\noindent
    {\bf #1 #2:}\if(#3\thmp\else\thmn#3\fi}
\def\neweqf#1$${\xdef #1{(\number\eqnumbo)}
    \eqno #1$$
    \global \advance \eqnumbo by 1}
\def\newthmf#1 #2: #3{    \global \advance \thmno by 1\xdef #2{\number\thmno}
    \medbreak\noindent
    {\bf #1 #2:}\if(#3\thmp\else\thmn#3\fi}
\def\forclose#1{\hfil\llap{$#1$}\hfilneg}
\def\newforclose#1{
	\ifsect\xdef #1{(\number\parano.\number\eqnumbo)}\else
	\xdef #1{(\number\eqnumbo)}\fi
	\hfil\llap{$#1$}\hfilneg
	\global \advance \eqnumbo by 1
	\iffinal\else\rsimb#1\fi}
\def\forevery#1#2$${\displaylines{\let\eqno=\forclose
        \let\neweq=\newforclose\hfilneg\rlap{$\qquad\quad\forall#1$}\hfil#2\cr}$$}
\def\noNota #1\par{}
\def\today{\ifcase\month\or
   January\or February\or March\or April\or May\or June\or July\or August\or
   September\or October\or November\or December\fi
   \space\number\year}
\def\inizia{\ifsect\let\neweq=\neweqt\else\let\neweq=\neweqf\fi
\ifsect\let\newthm=\newthmt\else\let\newthm=\newthmf\fi}
\def\bititolo{\empty}
\gdef\begin #1 #2\par{\xdef\titolo{#2}
\ifsect\let\neweq=\neweqt\else\let\neweq=\neweqf\fi
\ifsect\let\newthm=\newthmt\else\let\newthm=\newthmf\fi
\iffinal\let\Nota=\noNota\fi
\centerline{\titlefont\titolo}
\if\bititolo\empty\else\medskip\centerline{\titlefont\bititolo}
\xdef\titolo{\titolo\ \bititolo}\fi
\bigskip
\centerline{\bigfont
\autore \ifvoid\notaautore\else\footnote{${}^1$}{\unhbox\notaautore}\fi}
\bigskip\if\istituto!\centerline{\today}\else
\centerline{\istituto}
\centerline{\indirizzo}
\centerline{\email}
\medskip
\centerline{#1\ \anno}\fi
\bigskip\bigskip
\ifsect\else\global\thmno=1\global\eqnumbo=1\fi}
\def\anno{2014}
\def\raggedleft{\leftskip2cm plus1fill \spaceskip.3333em \xspaceskip.5em
\parindent=0pt\relax}
\font\titlefont=cmssbx10 scaled \magstep1
\font\bigfont=cmr12
\font\eightrm=cmr8
\font\sc=cmcsc10
\font\bbr=msbm10
\font\sbbr=msbm7
\font\ssbbr=msbm5
\font\msam=msam10

\font\bfm=cmmib10

\nopagenumbers
\binoppenalty=10000
\relpenalty=10000
\newfam\amsfam
\textfont\amsfam=\bbr \scriptfont\amsfam=\sbbr \scriptscriptfont\amsfam=\ssbbr
\newfam\boldifam
\textfont\boldifam=\bfm
\let\de=\partial
\let\eps=\varepsilon
\let\phe=\varphi
\def\Hol{\mathop{\rm Hol}\nolimits}

\def\Re{\mathop{\rm Re}\nolimits}
\def\Im{\mathop{\rm Im}\nolimits}

\def\id{\mathop{\rm id}\nolimits}

\def\bigoperp{\mathop{\hbox{$\bigcirc\kern-11.8pt\perp$}}\limits}
\def\Klim{\mathop{\hbox{$K$-$\lim$}}\limits}
\def\rKlim{\mathop{\hbox{$K'$-$\lim$}}\limits}

\mathchardef\void="083F
\mathchardef\ellb="0960
\mathchardef\taub="091C
\def\C{{\mathchoice{\hbox{\bbr C}}{\hbox{\bbr C}}{\hbox{\sbbr C}}
{\hbox{\sbbr C}}}}
\def\R{{\mathchoice{\hbox{\bbr R}}{\hbox{\bbr R}}{\hbox{\sbbr R}}
{\hbox{\sbbr R}}}}

\newcount\notitle
\notitle=1
\headline={\ifodd\pageno\rhead\else\lhead\fi}
\def\rhead{\ifnum\pageno=\notitle\iffinal\hfill\else\hfill\tt Version
\the\versiono; \the\day/\the\month/\the\year\fi\else\hfill\eightrm\titolo\hfill
\folio\fi}
\def\lhead{\ifnum\pageno=\notitle\hfill\else\eightrm\folio\hfill\autore\hfill
\fi}
\newbox\bibliobox
\def\setref #1{\setbox\bibliobox=\hbox{[#1]\enspace}
    \parindent=\wd\bibliobox}
\def\biblap#1{\noindent\hang\rlap{[#1]\enspace}\indent\ignorespaces}
\def\art#1 #2: #3! #4! #5 #6 #7-#8 \par{\biblap{#1}#2: {\sl #3\/}.
    #4 {\bf #5} (#6)\if.#7\else, \hbox{#7--#8}\fi.\par\smallskip}
\def\book#1 #2: #3! #4 \par{\biblap{#1}#2: {\bf #3.} #4.\par\smallskip}
\def\coll#1 #2: #3! #4! #5 \par{\biblap{#1}#2: {\sl #3\/}. In {\bf #4,}
#5.\par\smallskip}
\def\pre#1 #2: #3! #4! #5 \par{\biblap{#1}#2: {\sl #3\/}. #4, #5.\par\smallskip}
\def\raggedleft{\leftskip2cm plus1fill \spaceskip.3333em \xspaceskip.5em
\parindent=0pt\relax}
\def\Nota #1\par{\medbreak\begingroup\Bittersweet\raggedleft
#1\par\endgroup\Black
\ifdim\lastskip<\medskipamount \removelastskip\penalty 55\medskip\fi}
%

%
%\newcount\defno
\def\smallsect #1. #2\par{\bigbreak\noindent{\bf #1.}\enspace{\bf #2}\par
    \global\parano=#1\global\eqnumbo=1\global\thmno=0%\global\defno=0\global\remno=0
    \nobreak\smallskip\nobreak\noindent\message{#2}}
\def\newdef #1\par{\global \advance \thmno by 1
    \medbreak
{\bf Definition \the\parano.\the\thmno:}\enspace #1\par
\ifdim\lastskip<\medskipamount \removelastskip\penalty 55\medskip\fi}
%\finalfalse
%\versiono=12
\finaltrue

\magnification\magstephalf

%%%% Abbreviazioni
\let\no=\noindent

\let\me=\medskip
\let\sm=\smallskip

\def\rKlim{\mathop{\hbox{\rm $K'$-lim}}\limits}

\def\autore{Marco Abate${}^1$\footnote{}{\eightrm 2010 Mathematics Subject Classification: Primary 37L05; Secondary 32A40, 32H50, 20M20.\hfill\break\indent Keywords: infinitesimal generators, semigroups of holomorphic mappings, Julia-Wolff-Carath\'eodory theorem, boundary behaviour.}, Jasmin Raissy${}^2$\footnote{${}^*$}{\eightrm Partially supported by the FIRB2012 grant ``Differential Geometry and Geometric Function Theory''.}}
\def\indirizzo{\vbox{\hfill${}^2$Institut de Math\'ematiques de Toulouse; UMR5219,
Universit\'e de Toulouse; CNRS,\hfill\break\null\hfill
UPS IMT, F-31062 Toulouse Cedex 9, France.
E-mail: jraissy@math.univ-toulouse.fr\hfill\null}}
\def\istituto{\vbox{\hfill${}^1$Dipartimento di Matematica, Universit\`a
di Pisa,\hfill\break\null\hfill Largo Pontecorvo 5, 56127 Pisa,
Italy. E-mail: abate@dm.unipi.it\hfill\null\vskip5pt}}
\def\email{}

\begin {June} A Julia-Wolff-Carath\'eodory theorem for infinitesimal generators in the unit ball

{{\narrower{\sc Abstract.}  
We prove a Julia-Wolff-Carath\'edory theorem on angular derivatives of infinitesimal generators of one-parameter semigroups of holomorphic self-maps of the unit ball $B^n\subset\C^n$, starting from results recently obtained by Bracci and Shoikhet. \par
}}

\smallsect 0. Introduction

The classical Fatou theorem says that a bounded holomorphic function~$f$ defined
on the unit disk~$\Delta\subset\C$ admits non-tangential limit at almost every
point of~$\de\Delta$, but it does not say anything about the behavior
of~$f(\zeta)$ as~$\zeta$ approaches a specific point~$\sigma$ of the boundary.
Of course, to be able to say something in this case one needs some hypotheses
on~$f$. For instance, one can assume that, in a very weak sense, $f(\zeta)$
approaches the boundary of~$\Delta$ at least as fast as~$\zeta$. It turns out that under this condition, not only~$f$, but even its derivative admits non-tangential limit. This
is the content of the classical {\it Julia-Wolff-Carath\'eodory theorem:}

\newthm Theorem \zJulia: (Julia-Wolff-Carath\'eodory)
Let $f\colon\Delta\to\Delta$ be a
bounded holomorphic function such that
$$
\liminf_{\zeta\to\sigma}{1-|f(\zeta)|\over1-|\zeta|}=\alpha<+\infty
\neweq\eqzalp
$$
for some $\sigma\in\de\Delta$. Then $f$ has non-tangential limit $\tau\in\de
\Delta$ at~$\sigma$, for all $\zeta\in\Delta$ one has 
$$
{|\tau-f(\zeta)|^2\over1-|f(\zeta)|^2}\le\alpha\,{|\sigma-\zeta|^2\over1-
	|\zeta|^2}\;,
\neweq\eqzJ
$$
and furthermore both the incremental ratio $\bigl(\tau-f(\zeta)
\bigr)\big/(\sigma-\zeta)$ and the derivative~$f'(\zeta)$ have 
non-tangential limit~$\alpha\bar\sigma\tau$ at~$\sigma$.

This results from the work of
several authors: Julia~[Ju1, Ju2], Wolff~[Wo], Carath\'eodory~[C], Landau and
Valiron~[L-V], R.~Nevanlinna~[N] and others (see, e.g., [B] and~[A1] for proofs,
history and applications).

As already noticed by Kor\'anyi and Stein ([Ko], [K-S], [St]) when they extended Fatou's theorem to several complex variables, for domains in $\C^n$ the notion of non-tangential limit is not the right one to consider. Actually, it turns out that for generalizing the Julia-Wolff-Carath\'eodory theorem from the unit disk to the unit ball $B^n\subset\C^n$ one needs two different notions of limit at the boundary, both stronger than non-tangential limit.

%To describe the right notion (more precisely, notions: we shall need two of them) it is worthwhile to take a closer look
%to the notion of non-tangential limit in~$\Delta$. First of all, the 
%non-tangential limit can be defined in two equivalent ways. We can say that a
%function $f\colon\Delta\to\C$ has {\it non-tangential limit}~$L\in\C$
%at~$\sigma\in\de\Delta$ if~$f\bigl(\gamma(t)\bigr)\to L$ as~$t\to 1^-$ for
%every curve $\gamma\colon[0,1)\to\Delta$ such that $\gamma(t)\to\sigma$
%non-tangentially as~$t\to1^-$. Or, we can say that $f$ has non-tangential
%limit~$L\in\C$ if $f(\zeta)\to L$ as~$\zeta\to\sigma$ staying inside any {\it
%Stolz region}~$K(\sigma, M)$ of {\it vertex}~$\sigma$ and~{\it
%amplitude}~$M>1$, where
%$$
%K(\sigma,M)=\left\{\zeta\in\Delta\biggm|{|\sigma-\zeta|\over1-|\zeta|}<M
%	\right\}\;.
%$$
%Since Stolz regions are angle-shaped nearby the vertex~$\sigma$, and the angle
%is going to~$\pi$ as~$M\to+\infty$, it is clear that these two definitions are
%equivalent; but, as we shall momentarily see, in $\C^n$ they give rise to different notions.

\def\autore{Marco Abate, Jasmin Raissy}
A function $f\colon B^n\to\C$ has non-tangential limit $L\in\C$ at a boundary point $p\in\de B^n$ 
if $f(z)\to L$ as $z\to p$ staying inside cones with vertex at~$p$; a stronger notion of limit can be obtained by using approach regions larger than cones. 

In the unit disk, as approach regions for the non-tangential limit one can use Stolz regions, since they are angle-shaped nearby the vertex. 
In the unit ball $B^n\subset\C^n$ the natural generalization of a Stolz region is the  {\it Kor\'anyi region} $K(p, M)$ of {\it vertex}~$p\in\de B^n$ and~{\it
amplitude}~$M>1$ given by
$$
K(p,M)=\left\{z\in B^n\biggm|{|1-\langle z,p\rangle|\over1-\|z\|}<M
	\right\}\;,
$$
where $\|\cdot\|$ denotes the euclidean norm and $\langle\cdot\,,\cdot\rangle$ the canonical hermitian product. We shall say that a function $f\colon B^n\to\C$ has {\it $K$-limit} (or {\it admissible limit})~$L\in\C$ at~$p\in\de B^n$,
and we shall write
$\Klim\limits_{z\to p} f(z)=L$,
if $f(z)\to L$ as~$z\to p$ staying inside any Kor\'anyi region~$K(\sigma, M)$. Since a Kor\'anyi region $K(p,M)$ approaches the boundary non-tangentially along the normal direction at~$p$ but tangentially along the complex tangential directions at~$p$, it turns out that having $K$-limit is stronger than having 
non-tangential limit. However, the best generalization of Julia's lemma to $B^n$ is the following result (proved by Herv\'e [H] in terms of non-tangential limits and by Rudin [R] in general):

\newthm Theorem \zJdue: Let $f\colon B^n\to B^m$ be a holomorphic map such that 
$$
\liminf_{z\to p}{1-\|f(z)\|\over1-\|z\|}=\alpha<+\infty\;,
$$
for some $p\in\de B^n$. Then $f$ admits $K$-limit $q\in\de B^m$ at~$p$, and furthermore for all $z\in B^n$ one has
$$
{|1-\langle f(z),q\rangle|^2\over 1-\|f(z)\|^2}\le \alpha {|1-\langle z,p\rangle|^2\over 1-\|z\|^2}\;.
$$

To obtain a complete generalization of the 
Julia-Wolff-Carath\'eodory for~$B^n$ one needs a different notion of limit, still stronger than non-tangential limit, but weaker than $K$-limit. 

A crucial one-variable result relating limits along curves and non-tangential limits is {\it Lindel\"of's theorem.} Given $\sigma\in\de\Delta$, a {\it $\sigma$-curve} is a continuous curve $\gamma\colon[0,1)\to\Delta$ such that $\gamma(t)\to\sigma$ as $t\to 1^-$. Then Lindel\"of [Li] proved that if a bounded holomorphic function $f\colon\Delta\to\C$ admits limit~$L\in\C$ along a given $\sigma$-curve then it admits limit $L$ along all non-tangential $\sigma$-curves --- and thus it has non-tangential limit~$L$ at~$\sigma$. 

Trying to generalize this theorem to several complex variables, \v Cirka [\v C] realized that for a bounded holomorphic function the existence of the limit along a (suitable) $p$-curve (where $p\in\de B^n$) implies not only the existence of the non-tangential limit, but also the existence of the limit along any curve belonging to a larger class of curves, including some tangential ones --- but it does not in general imply the existence of the $K$-limit. To describe the version (due to Rudin [R]) of \v Cirka's result we shall need in this paper, let us introduce a bit of terminology. 

Let $p\in\de B^n$. As before, a {\it $p$-curve} is a continuous curve $\gamma\colon[0,1)\to B^n$ such that
$\gamma(t)\to p$ as $t\to 1^-$. A $p$-curve is {\it special} if
$$
\lim_{t\to 1^-}{\|\gamma(t)-\langle\gamma(t),p\rangle p\|^2\over 1-|\langle\gamma(t),p\rangle|^2}=0\;;
\neweq\eqzspec
$$
and, given $M>1$, it is {\it $M$-restricted} if
$$
{|1-\langle\gamma(t),p\rangle|\over 1-|\langle\gamma(t),p\rangle|}<M
$$
for all $t\in[0,1)$. We also say that $\gamma$ is {\it restricted} if it is $M$-restricted for some~$M>1$. In other words, $\gamma$ is restricted if and only if $t\mapsto\langle\gamma(t),p\rangle$ goes to~1 non-tangentially in~$\Delta$. 

It is not difficult to see that non-tangential curves are special and restricted; on the other hand,
a special restricted curve approaches the boundary non-tangentially along the normal direction, but it can approach the boundary tangentially along complex tangential directions. Furthermore,
a special $M$-restricted $p$-curve is eventually contained in any $K(p,M')$ with $M'>M$,
and conversely a special $p$-curve eventually contained in $K(p,M)$ is $M$-restricted. 
However, $K(p,M)$ can contain $p$-curves that are restricted but not special: for these curves the limit in \eqzspec\ might be a strictly positive number. 

With these definitions in place, we shall say that a function $f\colon B^n\to\C$ has
{\it restricted $K$-limit} (or {\it hypoadmissible limit})~$L\in\C$ at~$p\in\de B^n$, and we shall write
$\rKlim_{z\to p} f(z)=L$,
if $f\bigl(\gamma(t)\bigr)\to L$ as~$t\to 1^-$ for any special restricted $p$-curve $\gamma\colon[0,1)\to B^n$. It is clear that the existence of the $K$-limit implies the existence of the restricted $K$-limit, that in turns implies the existence of the non-tangential limit; but none of these implications can be reversed (see, e.g., [R] for examples in the ball). 

Finally, we say that a function $f\colon B^n\to\C$ is
{\it $K$-bounded} at~$p\in\de B^n$ if it is bounded in any Kor\'anyi region $K(p,M)$, where the bound can depend on~$M>1$. Then the version of \v Cirka's generalization of Lindel\"of's theorem we shall need is the following:

\newthm Theorem \zCR: (Rudin [R])
Let $f\colon B^n\to\C$ be a holomorphic function $K$-bounded at $p\in\de B^n$. Assume there is a special restricted $p$-curve $\gamma^o\colon[0,1)\to B^n$ such that $f\bigl(\gamma^o(t)\bigr)\to L\in\C$ as~$t\to 1^-$. Then $f$ has restricted $K$-limit $L$ at~$p$.

We can now deal with the generalization of the Julia-Wolff-Carath\'eodory theorem to several complex variables. With respect to the one-dimensional case there is an obvious difference:
instead of only one derivative we have to consider a whole (Jacobian) matrix of
them, and there is no reason they should all behave in the same way. And
indeed they do not, as shown in
Rudin's version of the Julia-Wolff-Carath\'eodory theorem for the unit ball:

\newthm Theorem \zJWCR: (Rudin [R])
Let $f\colon B^n\to B^m$ be a holomorphic map such that 
$$
\liminf_{z\to p}{1-\|f(z)\|\over1-\|z\|}=\alpha<+\infty\;,
$$
for some $p\in\de B^n$. Then $f$ admits 
$K$-limit $q\in\de B^m$ at $p$. Furthermore, if we set
$f_q(z)=\bigl\langle f(z),p
\bigr\rangle q$ and
denote by~$df_z$ the differential of~$f$ at~$z$, we have:
\item{\rm (i)} the function
$\bigl(1-\bigl\langle f(z),q\bigr\rangle\bigr)\big/(1-\langle z,p\rangle)$ is $K$-bounded and has restricted $K$-limit~$\alpha$ at~$p$;
\item{\rm (ii)} the map $(f(z)-f_q(z))/(1-\langle z,p\rangle)^{1/2}$ is $K$-bounded and has restricted
$K$-limit~$O$ at~$p$;
\item{\rm (iii)} the function $\bigl\langle df_z(p),q\bigr\rangle$ is $K$-bounded and has restricted
$K$-limit~$\alpha$ at~$p$;
\item{\rm (iv)} the map $(1-\langle z,p\rangle)^{1/2}d(f-f_q)_z(p)$ is $K$-bounded and has restricted
$K$-limit~$O$ at~$p$;
\item{\rm (v)} if $v$ is any vector orthogonal to~$p$, the function
$\bigl\langle df_z (v),q\bigr\rangle\big/(1-\langle z,p\rangle)^{1/2}$ 
is $K$-bounded and has restricted $K$-limit~$0$ at~$p$;
\item{\rm (vi)} if $v$ is any vector orthogonal to~$p$, the map
$d(f-f_q)_z (v)$ is $K$-bounded at~$p$.

In the last twenty years this theorem (as well as Theorems~\zJdue\ and~\zCR) has been extended to domains much more general than the unit ball: for instance, strongly pseudoconvex domains, convex domains of finite type, and polydisks (see, e.g., [A1], [A2], [A3], [A5], [AT], [A6], [AMY] and references therein). But in this paper we are interested in a different kind of generalization,
that we are now going to describe.

Let $\Hol(B^n,B^n)$ denote the space of holomorphic self-maps of~$B^n$, endowed with the usual compact-open topology. A {\it one-parameter semigroup} of holomorphic self-maps of~$B^n$ is a continuous semigroup homomorphism $\Phi\colon\R^+\to\Hol(B^n,B^n)$. In other words, writing $\phe_t$ instead of~$\Phi(t)$, we have $\phe_0=\id_{B^n}$, the map $t\mapsto\phe_t$ is continuous, and the semigroup property
$\phe_t\circ\phe_s=\phe_{t+s}$
holds (see, e.g., [A1], [RS2] or [S] for an introduction to the theory of one-parameter semigroups of holomorphic maps). 

One-parameter semigroups can be seen as the flow of a vector field (see, e.g., [A4]). Indeed, given a 
one-parameter semigroup $\Phi$, it is possible to prove that it exists a holomorphic map 
$G\colon B^n\to\C^n$, the {\it infinitesimal generator} of the semigroup, such that
$$
{\de\Phi\over\de t}=G\circ\Phi\;.
\neweq\eqzinfgen
$$
The infinitesimal generator can be obtained by the following formula:
$$
G(z)=\lim_{t\to 0^+}{\phe_t(z)-z\over t}\;.
\neweq\eqzigb
$$

\newrem In some papers (e.g., in [ERS] and [RS1]), the infinitesimal generator is defined as the solution of
the equation
$$
{\de\Phi\over\de t}+G\circ\Phi=O\;,
$$
that is with a change of sign with respect to our definition. This should be kept in mind when reading the literature on this subject.

Somewhat surprisingly, in 2008 Elin, Reich and Shoikhet [ERS] discovered a Julia's lemma for infinitesimal generators, just assuming that the radial limit of the generator at a point $p\in\de B^n$ vanishes (roughly speaking, this means that $p$ is a boundary fixed point for the associated semigroup):

\newthm Theorem \ERS: ([ERS, Theorem p. 403]) Let $G\colon B^n\to\C^n$ be the infinitesimal generator on $B^n$ of the one-parameter semigroup $\Phi=\{\phe_t\}$, and let $p\in\de B^n$ be such that
$$
\lim_{t\to 1^-}G(tp)=O\;.
\neweq\eqzbrnp
$$
Then the following assertions are equivalent:
\smallskip
\itm{(I)} we have
%$$
%\limsup_{t\to 1^-} \Re{\langle G(tp),p\rangle\over 1-t}>-\infty
%$$
%or, equivalently,
$$
\alpha=\liminf_{t\to 1^-} \Re{\langle G(tp),p\rangle\over t-1}<+\infty\;;
$$
%\itm{(II)} the limit
%$$
%\alpha:=\lim_{t\to 1^-} \Re{\langle G(tp),p\rangle\over t-1}
%$$
%exists and it is finite;
\itm{(II)} we have
$$
\beta= 2\sup_{z\in B^n}\Re \left[{\langle G(z),z\rangle\over 1-\|z\|^2}-{\langle G(z),p\rangle\over 1-\langle z,p\rangle}
\right]<+\infty\;;
$$
\itm{(III)} there exists $\gamma\in\R$ such that for all $z\in B^n$ we have
$$
{|1-\langle\phe_t(z),p\rangle|^2\over 1-\|\phe_t(z)\|^2}\le e^{\gamma t}
{|1-\langle z,p\rangle|^2\over 1-\|z\|^2}\;.
$$
%or, equivalently,
%$$
%u_p\bigl(\phi_t(z)\bigr)\le e^{-\gamma t} u_p(z)\;,
%$$
%where 
%$$
%d_p(z)={|1-\langle z,p\rangle|^2\over 1-\|z\|^2}=-{1\over u_p(z)}\;.
%$$
\smallskip
\noindent Furthermore, if any of these assertions holds then $\alpha=\beta =\inf\gamma$ and
we also have
$$
\lim_{t\to 1^-}{\langle G(tp),p\rangle\over t-1}=\beta\;.
\neweq\eqzbeta
$$
%\itm{(iii)} for all $z$,~$w\in B^n$ it holds
%$$
%\Re\left[{\langle G(z),z\rangle\over 1-\|z\|^2}+{\langle w, G(w)\rangle\over 1-\|w\|^2}\right]
%\le \Re\left[{\langle G(z),w\rangle+\langle z,G(w)\rangle\over 1-\langle z,w\rangle}\right]\;.
%$$

If \eqzbrnp\ and any (and hence all) of the equivalent conditions (I)--(III) holds we say that~$p\in\de B^n$ is a {\it boundary regular null point} of~$G$ with {\it dilation}~$\beta\in\R$. 

This result strongly suggests that one should try and prove a Julia-Wolff-Carath\'eodory theorem 
for infinitesimal generators along the line of Rudin's Theorem~\zJWCR. This has been partially achieved by Bracci and Shoikhet [BS], who proved the following 

\newthm Theorem \zJWCBS: ([BS]) Let $G\colon B^n\to\C^n$ be an infinitesimal generator on 
$B^n$ of a one-parameter semigroup, and let $p\in\de B^n$. Assume that
$$
{\langle G(z),p\rangle\over \langle z,p\rangle-1}\hbox{\qquad is $K$-bounded at $p$}
\neweq\eqzassdue
$$
and
$$
{G(z)-\langle G(z),p\rangle p\over (\langle z,p\rangle-1)^{1/2}}\hbox{\qquad is $K$-bounded at $p$.}
\neweq\eqzassdueb
$$
Then $p$ is a boundary regular null point for $G$. Furthermore, if $\beta$ is the dilation of $G$ at~$p$ then
\smallskip
\itm{(i)} the function
$\langle G(z),p\rangle\big/(\langle z,p\rangle-1)$ (is $K$-bounded and) has restricted $K$-limit~$\beta$ at~$p$;
\itm{(ii)} if $v$ is a vector orthogonal to $p$, the function $\langle G(z),v\rangle/(\langle z,p\rangle-1)^{1/2}$ is $K$-bounded at~$p$;
\itm{(iii)} the function $\langle dG_z(p),p\rangle$ is $K$-bounded and has restricted $K$-limit $\beta$ at~$p$;
\itm{(iv)} if $v$ is a vector orthogonal to~$p$, the function $(\langle z,p\rangle-1)^{1/2}\langle dG_z(p),v\rangle$ is $K$-bounded at~$p$;
\itm{(v)} if $v$ is a vector orthogonal to~$p$, the function $\langle dG_z(v),p\rangle\bigm/ (\langle z,p\rangle-1)^{1/2}$ is $K$-bounded at $p$.
\itm{(vi)} if $v_1$ and $v_2$ are vectors orthogonal to~$p$ the function $\langle dG_z(v_1),v_2\rangle$ is $K$-bounded at~$p$.

\newrem In the context of holomorphic maps, conditions \eqzassdue\ and \eqzassdueb\ are a consequence of (the equivalent of) condition (I) in Theorem~\ERS,
and indeed they appear as part of Theorem~\zJWCR.(i) and (ii); however, the proof in that setting uses in an essential way the fact that there we are dealing with holomorphic {\sl self-maps} of the ball. 
On the other hand, in our context, \eqzassdueb\ is {\sl not} a consequence of Theorem~\ERS.(I), as Example~1.2 shows, and \eqzassdue\ too seems to be stronger than Theorem~\ERS.(I); see also similar comments in [BS, Section~4.1]. Thus we have to assume \eqzassdue\ and \eqzassdueb\ as separate hypotheses. Furthermore, Example~1.2 also shows that the exponent $1/2$ might not necessarily be the right one to consider in the setting of infinitesimal generators. 

\newrem The assertions in Theorem~\zJWCBS.(i), (iii) and (v) follow just assuming \eqzassdue\ and that $G(tp)\to O$ as~$t\to 1^-$ (see [BS, Proposition~4.1]).

\newrem The assertions in Theorem~\zJWCBS\ (and in Theorem~0.12 below) have been numbered so as to reflect the similarities with the assertions in Theorem~\zJWCR. To see this, first of all notice that a boundary regular null point of~$G$ is a boundary fixed point of the associated semigroup~$\{\phe_t\}$. So in any comparison we must take ($m=n$ and) $q=p$ in Theorem~\zJWCR; in particular, the analogies between assertions (iii) and (v) in the two statements are obvious. Furthermore we can write
$$
{1-\langle \phe_t(z),p\rangle\over 1-\langle z, p\rangle}={\bigl\langle \phe_t(z)-p,p\bigr\rangle\over\langle z,p\rangle-1}={\bigl\langle \phe_t(z)-z,p\bigr\rangle\over\langle z,p\rangle-1}+1\;,
$$
and thus recalling \eqzigb\ it is clear that Theorem~\zJWCBS.(i) is the analogue of Theorem~\zJWCR.(i). Moreover, if $\{v_2,\ldots,v_n\}$ is an orthornormal basis of the vector space orthogonal to~$p$ we can write 
$$
G(z)-\bigl\langle G(z),p\bigr\rangle p =\sum_{j=2}^n \bigl\langle G(z),v_j\bigr\rangle v_j\;;
$$
therefore 
$$
d(G-\langle G,p\rangle p)_z(\cdot)=\sum_{j=2}^n\bigl\langle dG_z(\cdot),v_j\bigr\rangle v_j
$$
and the analogies between Theorem~\zJWCR.(ii), (iv) and~(vi) and the corresponding statements in Theorem~\zJWCBS\ become evident.

What is missing in Theorem~\zJWCBS\ to obtain a complete analogue of Theorem~\zJWCR\ is statements about restricted $K$-limits in cases (ii), (iv) and (v); the aim of this paper is exactly to provide those statements. It turns out that there is an obstruction, parallel to the one telling apart $K$-limits and restricted $K$-limits: as better described in Section~1, the curves one would like to use for obtaining the exponent $1/2$ in the statements are restricted but not special, in the sense that the limit in~\eqzspec\ is a strictly positive (though finite) number. We are thus led to consider exponents $\gamma <1/2$: this is not just a technical problem, but an inevitable feature of the theory, and in this way we actually widen the applicability of our results, as Example~1.2 shows. 

Our first main theorem then is:

\newthm Theorem \zJWCAR: Let $G\colon B^n\to\C^n$ be an infinitesimal generator on $B^n$ of a one-parameter semigroup, and let $p\in\de B^n$. Assume that
$$
{\langle G(z), p\rangle\over \langle z,p\rangle-1}\quad\hbox{and}\quad
{G(z)-\langle G(z),p\rangle p\over(\langle z,p\rangle-1)^{\gamma}} 
$$
are $K$-bounded at $p$ for some $0<\gamma<1/2$. Then $p\in\de B^n$ is a boundary regular null point  for $G$. Furthermore, if $\beta$ is the dilation of $G$ at~$p$ then:
\itm{(i)} the function
$\langle G(z),p\rangle\big/(\langle z,p\rangle-1)$ (is $K$-bounded and) has restricted $K$-limit~$\beta$ at~$p$;
\itm{(ii)} if $v$ is a vector orthogonal to $p$, the function $\langle G(z),v\rangle/(\langle z,p\rangle-1)^{\gamma}$ is $K$-bounded and has restricted $K$-limit $0$ at~$p$;
\itm{(iii)} the function $\langle dG_z(p),p\rangle$ is $K$-bounded and has restricted $K$-limit $\beta$ at~$p$;
\itm{(iv)} if $v$ is a vector orthogonal to~$p$, the function $(\langle z,p\rangle-1)^{1-\gamma}\langle dG_z(p),v\rangle$ is $K$-bounded and has restricted $K$-limit~$0$ at~$p$;
\itm{(v)} if $v$ is a vector orthogonal to~$p$, the function $\langle dG_z(v),p\rangle\bigm/ (\langle z,p\rangle-1)^{\gamma}$ is $K$-bounded and has restricted $K$-limit $0$ at $p$.
\itm{(vi)} if $v_1$ and $v_2$ are vectors orthogonal to~$p$ the function $(\langle z,p\rangle-1)^{1/2-\gamma}\langle dG_z(v_1),v_2\rangle$ is $K$-bounded\break\indent at~$p$.

%\smallskip
%\itm{(ii)} the map
%$$
%{G(z)-\langle G(z),p\rangle p\over(1-\langle z,p\rangle)^{\gamma}}
%$$
%(is $K$-bounded and) has restricted $K$-limit~$O$ at~$p$;
%\itm{(iv)} if $v$ is a vector orthogonal to~$p$ the function $(1-\langle z, p\rangle)^{1-\gamma}{\langle dG_z(p), v\rangle}$ is $K$-bounded and has restricted $K$-limit~$0$ at~$p$;
%\itm{(v)} if $v$ is a vector orthogonal to~$p$, the function $\langle dG_z(v),p\rangle\bigm/ (1-\langle z,p\rangle)^\gamma$ is $K$-bounded at $p$ and has restricted $K$-limit~$0$ at~$p$;\Red
%\itm{(vi)} if $v_1$ and $v_2$ are vectors orthogonal to~$p$ the function $\langle dG_z(v_1),v_2\rangle$ is $K$-bounded.\Black
% and for every $v\perp p$ the quotients 
%$$
%{\langle dG_z(v), p\rangle\over (1-\langle z, p\rangle)^{\gamma}}\;,
%\leqno{\rm(i)}
%$$
%and
%$$
%(1-\langle z, p\rangle)^{1-\gamma}{\langle dG_z(p), v\rangle}
%\leqno{\rm(ii)}
%$$
%are $K$-bounded at~$p$. Moreover, if
%$$
%\rKlim_{z\to p}{G(z)-\langle G(z),p\rangle p\over(1-\langle z,p\rangle)^{\gamma}}=0\;,
%\neweq\eqintrouno
%$$
%then
%$$
%\rKlim_{z\to p}{\langle dG_z(v), p\rangle\over (1-\langle z, p\rangle)^{\gamma}}=0\;,
%\neweq\eqintrodue
%$$
%and
%$$
%\rKlim_{z\to p}\, (1-\langle z, p\rangle)^{1-\gamma}{\langle dG_z(p), v\rangle}=0\;.
%\neweq\eqintrotre
%$$
%In particular, \eqintrouno\ holds for $0\le \gamma<1/2$, and for $\gamma=1/2$ with $p$ a H\"older boundary null point.

An exact analogue of Theorem~\zJWCR\ would be with $\gamma=1/2$; we can obtain such a statement by assuming a slightly stronger hypothesis on the infinitesimal generator.
Under the assumptions of Theorem~\zJWCBS\ we know that
$$
{\bigl\langle G\bigl(\sigma(t)\bigr),p\bigr\rangle\over \langle \sigma(t),p\rangle-1} = \beta+o(1)
\neweq\eqJC
$$
as $t\to 1^-$ for any special restricted $p$-curve $\sigma\colon[0,1)\to B^n$.
Following ideas introduced in [ESY], [EKRS] and [EJ] in the context of the unit disk, we shall say that $p$ is a {\it H\"older boundary null point} if there is $\alpha>0$ such that 
$$
{\bigl\langle G\bigl(\sigma(t)\bigr),p\bigr\rangle\over \langle \sigma(t),p\rangle-1} = \beta+o\bigl((1-t)^\alpha\bigr)
\neweq\eqH
$$
for any special restricted $p$-curve $\sigma\colon[0,1)\to B^n$ such that $\langle\sigma(t),p\rangle\equiv t$. Then our second main theorem is:

\newthm Theorem \zJWCARb: Let $G\colon B^n\to\C^n$ be an infinitesimal generator on $B^n$ of a one-parameter semigroup, and let $p\in\de B^n$. Assume that
$$
{\langle G(z), p\rangle\over \langle z,p\rangle-1}\quad\hbox{and}\quad
{G(z)-\langle G(z),p\rangle p\over(\langle z,p\rangle-1)^{1/2}} 
$$
are $K$-bounded at $p$, and that $p$ is a H\"older boundary null point. Then the statement of Theorem~\zJWCAR\ holds with $\gamma=1/2$.

We end this paper giving examples of infinitesimal generators with a H\"older boundary null point and satisfying the hypotheses of Theorem~\zJWCARb.

\me\no{\it Acknowledgments.} We gratefully thank Filippo Bracci for several useful discussions about the construction of Example~1.2, and David Shoikhet for pointing out to us references [ESY], [EKRS] and [EJ].

\smallsect 1. Proofs

This section is devoted to the proofs of Theorem~\zJWCAR\ and Theorem~\zJWCARb.

\sm\noindent{\sl Proof of Theorem~\zJWCAR\/.\enspace} Our hypotheses ensure that $\lim_{t\to1^-} G(tp) = O$ and therefore, thanks to Theorem~\ERS\ we have that $p$ is a boundary regular null point for $G$. Let $\beta\in \R$ be the dilation of~$G$ at $p$.

\sm(i) This follows immediately from our hypotheses, thanks to Theorems~\zCR\ and~\ERS.

\sm (ii) Given a vector $v$ orthogonal to $p$, the $K$-boundedness of the function $\langle G(z),v\rangle/(\langle z,p\rangle-1)^{\gamma}$ follows immediately from that of $\bigl(G(z)-\langle G(z),p\rangle p\bigr)/(\langle z,p\rangle-1)^{\gamma}$. Analogously, to prove that the restricted $K$-limit at $p$ is zero, it suffices to prove  
$$
\rKlim_{z\to p}{G(z)-\langle G(z),p\rangle p\over(\langle z,p\rangle-1)^{\gamma}}=0\;.
\neweq\equnog
$$ 
Without loss of generality, we can assume $p=e_1$, and we write $z=(z_1,z')$ with $z'=(z_2,\ldots,z_n)$ for points in $\C^n$. In particular, we can replace $G(z)-\langle G(z),p\rangle p$ by $G(z)'= (G_2(z),\dots, G_n(z))$ in the statement we would like to prove, and by Theorem~\zCR\ to get the assertion it suffices to show that 
$$
\lim_{t\to 1^-}{G_j(te_1)\over(t-1)^{\gamma}}=0
\neweq\equunoc
$$ 
for all $j=2,\dots,n$.

Since $G$ is an infinitesimal generator with boundary regular null point $e_1$ having dilation $\beta\in\R$, Theorem~\ERS\ implies that 
$$
\Re\left[{\bigl\langle G(z),z\bigr\rangle\over 1-\|z\|^2}-{G_1(z)\over 1-z_1}\right]\le {\beta\over 2}\;
\neweq\eqdue
$$
for any $z\in B^n$.

%We know that $G(z)'/(z_1-1)^{\gamma}$ is $K$-bounded by hypothesis, therefore, by \v Cirka's theorem (see, e.g., [R, Theorem 8.4.8]) to prove \equnog\ it suffices to show that 
%$$
%\lim_{t\to 1^-}{G_j(te_1)\over(t-1)^{\gamma}}=0
%$$ 
%for all $j=2,\dots,n$.
%
Given $j\in\{2,\ldots,n\}$, fix $0<\eps<1$ and $\theta\in\R$; for $t\in(0,1)$, set 
$$
z_t = te_1 + e^{-i\theta} \eps(1-t)^{1-\gamma}e_j\in B^n\;.
$$
In particular, $t\mapsto z_t$ is a special restricted $e_1$-curve, and we have 
$$
1-\|z_t\|^2= (1-t)(1+t-\eps^2(1-t)^{1-2\gamma})\;.
$$

Now, \eqdue\ evaluated in $z_t$ becomes
$$
\Re\left[{ t G_1(z_t) + e^{i\theta}\eps (1-t)^{1-\gamma}G_j(z_t) \over 1-\|z_t\|^2}-{G_1(z_t)\over 1-\langle z_t, e_1\rangle}\right]\le {\beta\over 2}\;.
$$
Therefore
$$
\eqalign{
\Re\left[{e^{i\theta}\eps (1-t)^{1-\gamma}G_j(z_t) \over 1-\|z_t\|^2}\right]
&\le {\beta\over 2} + \Re\left[{G_1(z_t)\over1-\langle z_t, e_1\rangle}\right] - t\Re\left[{G_1(z_t)\over 1-\|z_t\|^2}\right]\cr
&= {\beta\over 2} + \Re\left[{G_1(z_t)\over1-\langle z_t, e_1\rangle}\right]\left(1 - {t (1-\langle z_t, e_1\rangle)\over 1-\|z_t\|^2}\right)\cr
%&= {\beta\over 2} + (-\beta+ o(1) )\left(1 - {t (1-\langle z_t, e_1\rangle)\over 1-\|z_t\|^2}\right)\cr
&= {\beta\over 2} + \Re\left[{G_1(z_t)\over1-\langle z_t, e_1\rangle}\right]\left(1 - {t \over 1+t -\eps^2(1-t)^{1-2\gamma}}\right)
\;.}
$$
Furthermore
$$
\eqalign{
\Re\left[{e^{i\theta} \eps (1-t)^{1-\gamma}G_j(z_t) \over 1-\|z_t\|^2}\right]
&={\eps(1-t)^{1-\gamma}(1-\langle z_t, e_1\rangle)^{\gamma}\over 1-\|z_t\|^2}{\Re[e^{i\theta}G_j(z_t)] \over (1-\langle z_t, e_1\rangle)^{\gamma}}\cr
%&= {\eps(1-t)\over 1-t^2 -\eps^2 (1-t)}{|G_j(z_t)| \over (1-\langle z_t, e_1\rangle)^{1/2}}
%\cr
&= {\eps\over 1+t -\eps^2(1-t)^{1-2\gamma}}{\Re[e^{i\theta}G_j(z_t)] \over (1-\langle z_t, e_1\rangle)^{\gamma}}
\;.\cr}
$$ 
Recalling Theorem~\zJWCBS, and in particular \eqJC, we get
$$
\eqalign{
{\Re[e^{i\theta}G_j(z_t)] \over (1-\langle z_t, e_1\rangle)^{\gamma}}
&\le \left( {\beta\over 2} + \Re\left[{G_1(z_t)\over1-\langle z_t, e_1\rangle}\right]\left(1 - {t \over 1+t -\eps^2(1-t)^{1-2\gamma}}\right)
\right){ 1+t -\eps^2(1-t)^{1-2\gamma}\over \eps}\cr
&={\beta\over 2}\cdot { 1+t -\eps^2(1-t)^{1-2\gamma}\over \eps} + \Re\left[{G_1(z_t)\over1-\langle z_t, e_1\rangle}\right]\left({ 1+t -\eps^2(1-t)^{1-2\gamma}\over \eps}  - {t \over \eps}\right)\cr
&={\beta\over 2}\cdot { 1+t -\eps^2(1-t)^{1-2\gamma}\over \eps}+\bigl(-\beta+o(1)\bigr)
\left({ 1-\eps^2(1-t)^{1-2\gamma}\over \eps}\right)\cr
&={\beta\over 2}{\eps^2(1-t)^{1-2\gamma}+t-1\over\eps}+o(1)\;.
\cr
%&={\beta\over 2}{\eps^2(1-t)^{\alpha-\gamma+1/2}-(1-t)^{3/2-\alpha-\gamma}\over\eps}+o\left({ 1-\eps^2(1-t)^{2\alpha}\over \eps(1-t)^{\alpha+\gamma-1/2}}\right)
%\cr
}
$$
%Taking $\gamma=1/2-\alpha$ we get
%$$
%{\Re[e^{i\theta}G_j(z_t)] \over (1-\langle z_t, e_1\rangle)^{\gamma}}\le
%{\beta\over 2}{\eps^2(1-t)^{2\alpha}-(1-t)\over\eps}+o(1)
%$$
Letting $t\to 1^-$ we obtain
$$
\limsup_{t\to 1^-}{\Re[e^{i\theta}G_j(z_t)] \over (1-\langle z_t, e_1\rangle)^{\gamma}}\le 0
$$
for all $\eps>0$ and $\theta\in\R$. Now letting $\eps\to 0^+$ we find
$$
\limsup_{t\to 1^-}{\Re[e^{i\theta}G_j(te_1)] \over (1-t)^{\gamma}}\le 0
$$
for all $\theta\in\R$, and this is possible if and only if
$$
\lim_{t\to 1^-}{G_j(te_1)\over (1-t)^{\gamma}}=0\;,
$$
and \equunoc\ follows. 
%that is
%$$
%\lim_{t\to 1^-}{G_j(te_1)\over (t-1)^{\gamma}}=0\;,
%$$
%as claimed.

\sm (iii) The proof %that $\langle dG_z(p),p\rangle$ is $K$-bounded at $p$ is an application of Cauchy's formula, and it 
is analogous to the one given in~[BS]; we recall it here for the sake of completeness. 

Without loss of generality, we can assume $p=e_1$. Let $M'>M> 1$ %and consider $K(e_1, M)$ the Kor\'anyi region of center $e_1$ and radius $M$. Let $M'>M$ 
and set $\delta:= {1\over 3} ({1\over M} - {1\over M'})$. Thanks to [R, Lemma 8.5.5], for any $z\in K(e_1, M)$ and $(\lambda, u')\in\C\times\C^{n-1}$ with $|\lambda|\le \delta |z_1-1|$ and $\|u'\|\le \delta |z_1-1|^{1/2}$, we have $(z_1+\lambda,z'+u')\in K(e_1, M')$. 

Now, fix $z\in K(e_1, M)$ and let $r = r(z) := \delta |z_1-1|$. By Cauchy's formula, we have
$$
\eqalign{
\langle dG_z(e_1),e_1\rangle
&=
{1\over 2\pi i}\int_{|\zeta| =r} {\langle G(z_1 +\zeta, z'), e_1\rangle \over\zeta^2}d\zeta\cr
&={1\over 2\pi}\int_{-\pi}^\pi {\langle G(z_1 +re^{i\theta}, z'), e_1\rangle \over
z_1 +r e^{i\theta} -1} {z_1 +r e^{i\theta}-1\over re^{i\theta}} d\theta\;.
}
$$
The first factor in the integral is bounded because $(z_1 +re^{i\theta}, z')\in K(e_1, M')$; furthermore, we also have $|({z_1 +r e^{i\theta}-1)/re^{i\theta}}|\le 1 + 1/\delta$, and thus we are done. 

To prove that the restricted $K$-limit at~$p$ is $\beta$, by  Theorem~\zCR\ it suffices to prove that 
$$
\lim_{t\to 1^-} \langle dG_{te_1}(e_1), e_1\rangle = \beta.
$$
Thanks to [BCD, Theorem 0.4], we have that $\lim_{t\to 1^-} {d \over dt}\left(G_1(te_1)\right)= \beta$, and then we are done, because ${d \over dt}\left(G_1(te_1)\right) = \langle dG_{te_1}(e_1), e_1\rangle$.

\sm (iv) Without loss of generality we can assume $p=e_1$ and $v=e_2$, so that the quotient
we would like to study is
$$
(z_1-1)^{1-\gamma}{\de G_2\over\de z_1}(z)\;.
$$
The proof of the $K$-boundedness is again an application of the Cauchy formula. %, and it is similar to the one of Theorem~\zJWCBS.(iv) for the case $\gamma=1/2$. We write it here for general $0<\gamma< 1/2$ for the sake of completeness.
As before, %Let $M> 1$ and consider $K(e_1, M)$ the Kor\'anyi region of center $e_1$ and radius $M$. 
let $M'>M>1$ and set $\delta:= {1\over 3} ({1\over M} - {1\over M'})$. Thanks to [R, Lemma 8.5.5], for any $z\in K(e_1, M)$ and $(\lambda, u')\in\C\times\C^{n-1}$ with $|\lambda|\le \delta |z_1-1|$ and $\|u'\|\le \delta |z_1-1|^{1/2}$, we have $(z_1+\lambda,z'+u')\in K(e_1, M')$. 

Now, fix $z\in K(e_1, M)$ and let $r = r(z) := \delta |z_1-1|$. By Cauchy's formula, we have
$$
\eqalign{
|z_1-1|^{1-\gamma}{\de G_2\over\de z_1}(z)
&=
{|z_1-1|^{1-\gamma}\over 2\pi i}\int_{|\zeta| =r} {G_2(z_1 +\zeta, z')\over\zeta^2}d\zeta\cr
&={1\over 2\pi \delta}\int_{-\pi}^\pi {G_2(z_1 +r e^{i\theta}, z')\over
|z_1 +r e^{i\theta} -1|^{\gamma}} \left| {z_1 +r e^{i\theta}-1\over z_1-1}\right|^{\gamma} {|z_1-1|\over|z_1-1|e^{i\theta}} d\theta\;.
}
$$
The choice of $r$ ensures that $(z_1 +\zeta, z')\in K(e_1,M')$; thus the first factor in the integral is bounded, and, since an easy computation shows that ${|z_1 +r e^{i\theta}-1|\over |z_1-1|}\le 1+\delta$, we are done.

To prove that the restricted $K$-limit at~$p$ vanishes, thanks to Theorem~\zCR, % i.e. \cite[Theorem 2.7.13]{A}/ 
it suffices to show that
$$
\lim_{t\to 1^-} (t-1)^{1-\gamma}{\de G_2\over\de z_1}(te_1) = 0\;.
\neweq\equnoterred
$$
Indeed, choose $\eps \in(0,1)$, and for any $t\in (0,1)$, let $\sigma_t\colon\eps\Delta\to B^n$ be defined by
$$
\sigma_t(\zeta) = (t + \zeta(1-t))e_1\;. 
$$
Then $\sigma_t(0) = t e_1$ and $\sigma_t'(0) = (1-t) e_1$. Moreover, for any $\zeta\in\eps\Delta$ we have
$$
{|1-t-\zeta(1-t)|\over1-|t+ \zeta(1-t)|} = {(1-t)|1-\zeta| \over 1-|1 - (1-t)(1- \zeta)|} \le {1+\eps \over1-\eps}\;.
$$
Therefore
$\sigma_t(\overline {\eps\Delta}) \subset K(e_1, M)$ for all $M >{1+\eps\over 1-\eps}$. In particular, for all $\theta\in \R$, the $e_1$-curve $t\mapsto \sigma_t(\eps e^{i\theta})$ is special and $M$-restricted. Now, 
$$
\eqalign{
(t-1)^{1-\gamma}{\de G_2\over\de z_1}(te_1) 
& = {1\over 2\pi} \int_{-\pi}^{\pi}{G_2(t+\eps (1-t)e^{i\theta}, O') \over (t+\eps(1-t)e^{i\theta}-1)^\gamma}{(t+\eps(1-t)e^{i\theta}-1)^\gamma\over \eps(1-t) e^{i\theta}}(t-1)^{1-\gamma} d\theta\cr
&= {-1\over 2\pi} \int_{-\pi}^{\pi}{G_2(t+\eps(1-t)e^{i\theta}, O') \over (t+\eps(1-t)e^{i\theta}-1)^\gamma} {(1-\eps e^{i\theta})^\gamma\over\eps e^{i\theta}}d\theta\;.
}
$$
The second factor of the integrand is bounded, and the first factor converges punctually and boundedly to $0$ as $t\to 1$, thanks to (ii); therefore \equnoterred\ follows from the dominated convergence theorem.

\sm (v) Without loss of generality we can assume $p=e_1$ and $v=e_2$, so that the quotient
we would like to study is
$$
{1\over(z_1-1)^\gamma}{\de G_1\over\de z_2}(z)\;.
$$
%For $\gamma=1/2$, the $K$-boundedness of this quotient has already been proved in Theorem~\zJWCBS.(v). The same proof works analogously for $0< \gamma< 1/2$, and thus we omit it here.
The proof of the $K$-boundedness is yet another application of the Cauchy formula.
%, and it is similar to the one of Theorem~\zJWCBS.(v) for the case $\gamma=1/2$. We write it here for general $0<\gamma< 1/2$ for the sake of completeness.
Let $M'>M>1$; %and consider $K(e_1, M)$ the Kor\'anyi region of center $e_1$ and radius $M$. Let $M'>M$ and 
set $\delta:= {1\over 3} ({1\over M} - {1\over M'})$, and $r = r(z) := \delta |z_1-1|^{1-\gamma}$;
[R, Lemma 8.5.5] ensures that if $z\in K(e_1, M)$ then $z +r e^{i\theta} e_2 \in K(e_1,M')$
for all~$\theta\in\R$. Then
%
%Thanks to [R, Lemma 8.5.5], for any $z\in K(e_1, M)$ and $(\lambda, u')\in\C\times\C^{n-1}$ with $|\lambda|\le \delta |z_1-1|$ and $\|u'\|\le \delta |z_1-1|^{1/2}$, we have $(z_1+\lambda,z'+u')\in K(e_1, M')$. 
%
%Now, fix $z\in K(e_1, M)$ and let $r = r(z) := \delta |z_1-1|^{1-\gamma}$. By 
Cauchy's formula yields
$$
\eqalign{
{1\over|z_1-1|^{\gamma}}{\de G_1\over\de z_2}(z)
&=
{1\over 2\pi i|z_1-1|^\gamma}\int_{|\zeta| =r} {G_1(z +\zeta e_2)\over\zeta^2}d\zeta\cr
&={1\over 2\pi \delta}\int_{-\pi}^\pi {G_1(z +r e^{i\theta} e_2)\over
|z_1 -1| e^{i\theta}} d\theta\;,
}
$$
and the $K$-boundness follows. 

Now we would like to prove that the restricted $K$-limit at~$p$ vanishes. %it suffices to prove that 
%$$
%\lim_{t\to 1^-}{1\over(t-1)^\gamma}{\de G_1\over\de z_2}(te_1)=0\;.
%$$ 
Let $\Phi\colon B^2\to B^n$ be given by $\Phi(\zeta,\eta)=\zeta e_1+\eta e_2$, and put $H=\Xi\circ\Phi$, where
$$
\Xi(z)={\bigl\langle G(z),z\bigr\rangle\over 1-\|z\|^2}-{G_1(z)\over 1-z_1}\;.
$$
Hence
$$
H(\zeta,\eta)={G_1(\zeta,\eta,0,\ldots,0)\bar\zeta+G_2(\zeta,\eta,0,\ldots,0)\bar\eta\over 1-|\zeta|^2-|\eta|^2}
-{G_1(\zeta,\eta,0,\ldots,0)\over 1-\zeta}\;.
$$
Now we expand $H$ in power series with respect to $\eta$:
$$
H(\zeta,\eta)=H(\zeta,0)+{\de H\over\de\eta}(\zeta,0)\eta+{\de H\over\de\bar\eta}(\zeta,0)\bar\eta
+O(|\eta|^2)\;.
\neweq\uquat
$$
We have
$$
H(\zeta,0)=G_1(\zeta,O')\left[{\bar\zeta\over1-|\zeta|^2}-{1\over1-\zeta}\right]=
-G_1(\zeta,O'){1\over1-|\zeta|^2}{1-\bar\zeta\over1-\zeta}\;;
$$
$$
{\de H\over\de\eta}(\zeta,0)={\de G_1\over\de z_2}(\zeta,O')\left[{\bar\zeta\over1-|\zeta|^2}-{1\over1-\zeta}\right]=-{\de G_1\over\de z_2}(\zeta,O'){1\over1-|\zeta|^2}{1-\bar\zeta\over1-\zeta}\;;
$$
and
$$
{\de H\over\de\bar\eta}(\zeta,0)={G_2(\zeta,O')\over1-|\zeta|^2}\;.
$$
Recalling \eqdue\ we get
$$
\eqalign{
{\beta\over 2}&\ge\Re H(\zeta,\eta)=\Re\left[H(\zeta,0)+{\de H\over\de\eta}(\zeta,0)\eta+{\de H\over\de\bar\eta}(\zeta,0)\bar\eta+O(|\eta|^2)\right]\cr
&={1\over 1-|\zeta|^2}\Re\left[-\left(G_1(\zeta,O')+\eta{\de G_1\over\de z_2}(\zeta,O')\right)
{|1-\zeta|^2\over(1-\zeta)^2}+G_2(\zeta,O')\bar\eta+O\bigl((1-|\zeta|^2)|\eta|^2\bigr)\right]\;,
\cr}
$$
and thus
$$
-{\beta\over 2}{1-|\zeta|^2 \over |1-\zeta|^2} \le\Re\left[{G_1(\zeta,O')\over(1-\zeta)^2}+{\eta\over(1-\zeta)^2}{\de G_1\over\de z_2}(\zeta,O')
- {\bar\eta G_2(\zeta,O')\over |1-\zeta|^2}+O\left({1-|\zeta|^2\over|1-\zeta|^2}|\eta|^2\right)\right]\;.
\neweq\eqtreb
$$

Fix $c>0$ and for $t\in[0,1)$ put
$$
\zeta_t = t + ic(1-t)\;.
$$
In particular,
$$
1-\zeta_t=(1-t)(1-ic)\;,\quad
|1-\zeta_t|=(1-t)(1+c^2)^{1/2}\quad\hbox{and}\quad {1\over 1-\zeta_t}={1\over 1-t}{1+ic\over 1+c^2}\;.
$$
It is easy to check that $\zeta_t\in \Delta$ if $1-t< 2/(1+c^2)$, and in this case
$$
1-|\zeta_t|^2=1-t^2-c^2(1-t)^2=(1-t)\bigl(1+t-(1-t)c^2\bigr)< 2(1-t)\;.
$$
Moreover, if
$1-t<1/(1+c^2)$ we have $1-|\zeta_t|^2> 1-t$, and thus we can find $\eta_t\in\C$ such that
$$
2(1-t)>1-|\zeta_t|^2>|\eta_t|^2>1-t\;;
$$
in particular, $(\zeta_t,\eta_t)\in B^2$, and we choose the argument of $\eta_t$ so that
$$
{\eta_t\over(1-\zeta_t)^2}{\de G_1\over\de z_2}(\zeta_t,O')=-\left|
{\eta_t\over(1-\zeta_t)^2}{\de G_1\over\de z_2}(\zeta_t,O')\right|\in\R^-\;.
$$
Now we compute \eqtreb\ in $(\zeta_t,\eta_t)$.
Multiplying by $|1-\zeta_t|^{2-\gamma}$ and dividing by $|\eta_t|$ we get
$$
\eqalign{
\left| {1\over(1-\zeta_t)^{\gamma}}{\de G_1\over\de z_2}(\zeta_t,O')\right|
&\le \Re\left[{G_1(\zeta_t,O')\over1-\zeta_t}{|1-\zeta_t|^{2-\gamma}\over(1-\zeta_t)|\eta_t|}\right]
+{|G_2(\zeta_t,O')|\over|1-\zeta_t|^{\gamma}}+O\left({1-|\zeta_t|^2\over|1-\zeta_t|^{\gamma}}|\eta_t|\right)\cr
&\quad+ {\beta\over 2}{1-|\zeta_t|^2 \over |1-\zeta_t|^\gamma |\eta_t|}\;.
\cr}
$$
Applying \eqJC\ we obtain
$$
\eqalign{
\biggl| {1\over(1-\zeta_t)^{\gamma}}&{\de G_1\over\de z_2}(\zeta_t,O')\biggr|\cr
&\le 
{|1-\zeta_t|^{2-\gamma}\over|\eta_t|}\Re\left[{-\beta+o(1)\over 1-t}{1+ic\over 1+c^2}\right]
+{|G_2(\zeta_t,O')|\over|1-\zeta_t|^{\gamma}}+O\left({1-|\zeta_t|^2\over|1-\zeta_t|^{\gamma}}|\eta_t|\right) \cr
&\quad+ {\beta\over 2}{1-|\zeta_t|^2 \over |1-\zeta_t|^{\gamma}|\eta_t|}\cr
&\le{(1-t)^{2-\gamma}(1+c^2)^{1-\gamma/2}\over(1-t)^{1/2}}{-\beta+o(1)\over(1-t)(1+c^2)}+
{|G_2(\zeta_t,O')|\over|1-\zeta_t|^{\gamma}}+O\left({1-|\zeta_t|^2\over|1-\zeta_t|^{\gamma}}|\eta_t|\right) \cr
&\quad+{\beta\over 2}{1-|\zeta_t|^2 \over (1-t)^{\gamma}(1+c^2)^{\gamma/2}|\eta_t|}\cr
&\le{(-\beta+o(1))(1-t)^{1/2-\gamma}\over(1+c^2)^{\gamma/2}}+{|G_2(\zeta_t,O')|\over|1-\zeta_t|^{\gamma}}+O\left({2(1-t)\over(1-t)^{\gamma}(1+c^2)^{\gamma/2}}\sqrt{2}(1-t)^{1/2}\right)\cr
&\quad+{\beta\over 2}{2(1-t) \over (1-t)^{\gamma}(1+c^2)^{\gamma/2}(1-t)^{1/2}}
\cr
&\le %{2\beta(1-t)^{1/2-\gamma}\over(1+c^2)^{\gamma/2}}+ 
o\bigl((1-t)^{1/2-\gamma}\bigr)+{|G_2(\zeta_t,O')|\over|1-\zeta_t|^{\gamma}}+O\bigl((1-t)^{3/2-\gamma}\bigr)\;.\cr}
$$
Since $t\mapsto\zeta_t e_1$ is a special restricted curve we can apply (ii) obtaining
$$
\limsup_{t\to 1}\biggl| {1\over(1-\zeta_t)^{\gamma}}{\de G_1\over\de z_2}(\zeta_t,O')\biggr|
\le 0\;.
%\cases{0& if $0<\gamma<1/2$,\cr
%\displaystyle {2\beta\over(1+c^2)^{1/4}}& if $\gamma=1/2$.}
$$
So %if $\gamma<1/2$ 
we get
$$
\lim_{t\to 1}{1\over(\zeta_t-1)^{\gamma}}{\de G_1\over\de z_2}(\zeta_t,O')=0
$$
and the assertion follows from Theorem~\zCR. 

\sm (vi) Without loss of generality we can assume $p=e_1$, $v_1=e_2$, and $v_2=e_3$, so that the function we would like to study is
$$
(z_1-1)^{{1\over 2}-\gamma}{\de G_3\over\de z_2}(z)\;.
$$
%For $\gamma=1/2$, the $K$-boundedness of this quotient has already been proved in Theorem~\zJWCBS.(v). The same proof works analogously for $0< \gamma< 1/2$, and thus we omit it here.
We argue as usual. 
%The proof of the $K$-boundedness is an application of the Cauchy formula, and it is similar to the one of Theorem~\zJWCBS.(vi) for the case $\gamma=1/2$. We write it here for general $0<\gamma< 1/2$ for the sake of completeness.
 
Let $M'>M>1$ %and consider $K(e_1, M)$ the Kor\'anyi region of center $e_1$ and radius $M$. Let $M'>M$ 
and set $\delta:= {1\over 3} ({1\over M} - {1\over M'})$. Thanks to [R, Lemma 8.5.5], for any $z\in K(e_1, M)$ and $u'\in\C^{n-1}$ with $\|u'\|\le \delta |z_1-1|^{1/2}$ we have $(z_1,z'+u')\in K(e_1, M')$. 

Now, fix $z\in K(e_1, M)$ and let $r = r(z) := \delta |z_1-1|^{1/2}$. By Cauchy's formula, we have
$$
\eqalign{
|z_1-1|^{{1\over 2}-\gamma}{\de G_3\over\de z_2}(z)
&=
{|z_1-1|^{{1\over 2}-\gamma}\over 2\pi i}\int_{|\zeta| =r} {G_3(z +\zeta e_2)\over\zeta^2}d\zeta\cr
&={1\over 2\pi \delta}\int_{-\pi}^\pi {G_3(z + r e^{i\theta} e_2)\over
|z_1 -1|^\gamma e^{i\theta}} d\theta\;.
}
$$
The choice of $r$ ensures that $z +r e^{i\theta} e_2 \in K(e_1,M')$, and the assertion follows from (ii).
\qedn

An accurate examination of the proof of the previous theorem reveals that the main point is the proof of part (ii). As soon as the statement of Theorem~\zJWCAR.(ii) holds for some $0<\gamma\le1/2$ (with $\gamma=1/2$ included) then the rest of the Theorem follows with the same $\gamma$ (again, $\gamma=1/2$ included). The proof of Theorem~\zJWCAR.(ii) we presented however breaks down for $\gamma=1/2$ because the curve 
$$
(0,1)\ni t\mapsto z_t = te_1 + e^{-i\theta} \eps(1-t)^{1-\gamma}e_j\in B^n
$$
is {\it not} special if $\gamma=1/2$; the limit \eqzspec\ is a strictly positive (though finite) number. 

\newrem Even assuming that the hypotheses of Theorem~\zJWCAR\ are satisfied with $\tilde\gamma \ge 1/2$, as explained above with this proof we can only obtain the thesis for all exponents $\gamma< 1/2$. 

Furthermore the exponent $1/2$, which is the natural one to consider in the setting of self-maps, it is not necessarily the right one for infinitesimal generators, as next example shows.

\newex Let $G\colon B^2\to \C^2$ be defined as 
$$
G(z,w) = (-z (1-z), -w (1-z)^{-\alpha})\;,
$$ 
with $0<\alpha<1/2$. It is easy to check that $G$ is an infinitesimal generator, since it vanishes at the origin and~$\Re\langle G(z,w), (z,w)\rangle \le 0$ for every $(z,w)\in B^2$. Moreover, $G$ satisfies the hypotheses of Theorem~\zJWCAR\ with $p=e_1$ and $\gamma = 1/2 -\alpha$, but $G_2(z,w)/(z-1)^\beta$ is not $K$-bounded for any $\beta> 1/2-\alpha$. Indeed, given $c \in(0,1)$,
all points of the form $(t, c\sqrt{1-t^2})$, with $t\in [0,1)$, belong to a Kor\'anyi region of vertex~$e_1$, whereas ${G_2(t,c\sqrt{1-t^2})/(t-1)^\beta}$ is not bounded as $t$ tends to $1$, for $1/2-\alpha -\beta < 0$. Furthermore $G_2(z,w)/(z-1)^\beta$ does not even have a restricted $K$-limit at $e_1$. In fact, choosing $\rho>1$ such that $\beta > \rho/2 -\alpha$, the curve $\sigma_\rho\colon [0,1)\to B^2$ defined by $\sigma_\rho(t) = (t, c(1-t^2)^{\rho/2})$, with $c\in (0,1)$, is a special restricted $e_1$-curve such that $G_2(\sigma_\rho(t))/(t-1)^\beta$ diverges as $t$ tends to $1$. 
This example can be easily generalized to any dimension.

%\newthm Corollary \corBS: Let $G\in\Hol(B^n,\C^n)$ be an infinitesimal generator, and let $p\in\de B^n$. If
%$$
%{\langle G(z), p\rangle\over \langle z,p\rangle-1}\quad\hbox{and}\quad {G(z)-\langle G(z),p\rangle p\over(\langle z,p\rangle-1)^{1/2}} 
%$$
%are $K$-bounded at $p$, then $p\in\de B^n$ is  a boundary regular null point, and 
%$$
%\rKlim_{z\to p}{G(z)-\langle G(z),p\rangle p\over(\langle z,p\rangle-1)^{\gamma}}=0\;,
%$$
%for any $0\le \gamma<1/2$.
%
%\pf Our hypotheses ensure that $\lim_{t\to1^-} G(tp) = O$ and therefore, Theorem~\zJWCBS\ gives us that $p$ is a boundary regular null point for $G$. Let $\beta\in \R$ be the dilation of $G$ at $p$.
%
%The rest of the statement is a straightforward consequence of the previous proposition, since if ${G(z)-\langle G(z),p\rangle p\over(\langle z,p\rangle-1)^{1/2}}$ is $K$-bounded at $p$, then ${G(z)-\langle G(z),p\rangle p\over(\langle z,p\rangle-1)^{\gamma}}$ is $K$-bounded at $p$ for any $0\le \gamma<1/2$.
%\qedn
%

On the other hand, we can get the statement with exponent $\gamma=1/2$ by using the notion of H\"older boundary null point, 
as follows:
 
\sm\noindent{\sl Proof of Theorem~\zJWCARb\/.\enspace} As explained above, it suffices to prove that 
$$
\lim_{t\to 1^-}{G_j(te_1)\over(t-1)^{1/2}}=0
\neweq\equunod
$$ 
for all $j=2,\dots,n$.

Let $\alpha>0$ be given by the definition of H\"older boundary null point; we can
clearly assume that $\alpha<1$.
%Since $G$ is an infinitesimal generator with boundary regular null point $e_1$ having dilation $\beta\in\R$, Theorem~\ERS\ implies that 
%$$
%\Re\left[{\bigl\langle G(z),z\bigr\rangle\over 1-\|z\|^2}-{G_1(z)\over 1-z_1}\right]\le {\beta\over 2}\;
%\neweq\eqdue
%$$
%for any $z\in B^n$.
%
%We know that $G(z)'/(z_1-1)^{\gamma}$ is $K$-bounded by hypothesis, therefore, by \v Cirka's theorem (see, e.g., [R, Theorem 8.4.8]) to prove \equnog\ it suffices to show that 
%$$
%\lim_{t\to 1^-}{G_j(te_1)\over(t-1)^{\gamma}}=0
%$$ 
%for all $j=2,\dots,n$.
%
Given $j\in\{2,\ldots,n\}$, fix $0<\eps<1$ and $\theta\in\R$; for $t\in(0,1)$, set 
$$
z_t = te_1 + e^{-i\theta} \eps(1-t)^{1/2+\alpha}e_j\in B^n\;.
$$
In particular, $t\mapsto z_t$ is a special restricted $e_1$-curve such that $\langle z_t,e_1\rangle\equiv t$, and we have 
$$
1-\|z_t\|^2= (1-t)(1+t-\eps^2(1-t)^{2\alpha})\;.
$$

Now, \eqdue\ evaluated in $z_t$ becomes
$$
\Re\left[{ t G_1(z_t) + e^{i\theta}\eps (1-t)^{1/2+\alpha}G_j(z_t) \over 1-\|z_t\|^2}-{G_1(z_t)\over 1-\langle z_t, e_1\rangle}\right]\le {\beta\over 2}\;.
$$
Therefore
$$
%\eqalign{
\Re\left[{e^{i\theta}\eps (1-t)^{1/2+\alpha}G_j(z_t) \over 1-\|z_t\|^2}\right]
%&\le {\beta\over 2} + \Re\left[{G_1(z_t)\over1-\langle z_t, e_1\rangle}\right] - t\Re\left[{G_1(z_t)\over 1-\|z_t\|^2}\right]\cr
%&= {\beta\over 2} + \Re\left[{G_1(z_t)\over1-\langle z_t, e_1\rangle}\right]\left(1 - {t (1-\langle z_t, e_1\rangle)\over 1-\|z_t\|^2}\right)\cr
%&= {\beta\over 2} + (-\beta+ o(1) )\left(1 - {t (1-\langle z_t, e_1\rangle)\over 1-\|z_t\|^2}\right)\cr
%&= 
\le {\beta\over 2} + \Re\left[{G_1(z_t)\over1-\langle z_t, e_1\rangle}\right]\left(1 - {t \over 1+t -\eps^2(1-t)^{2\alpha}}\right)\;.
%}
$$
Furthermore
$$
\eqalign{
\Re\left[{e^{i\theta} \eps (1-t)^{1/2+\alpha}G_j(z_t) \over 1-\|z_t\|^2}\right]
&={\eps(1-t)^{1/2+\alpha}(1-\langle z_t, e_1\rangle)^{1/2}\over 1-\|z_t\|^2}{\Re[e^{i\theta}G_j(z_t)] \over (1-\langle z_t, e_1\rangle)^{1/2}}\cr
&= {\eps(1-t)^{1+\alpha}\over 1-t^2 -\eps^2 (1-t)^{1+2\alpha}}{|G_j(z_t)| \over (1-\langle z_t, e_1\rangle)^{1/2}}
\cr
&
= {\eps(1-t)^\alpha\over 1+t -\eps^2(1-t)^{2\alpha}}{\Re[e^{i\theta}G_j(z_t)] \over (1-\langle z_t, e_1\rangle)^{1/2}}
\;.\cr}
$$ 
Using \eqH\ we then get
$$
\eqalign{
{\Re[e^{i\theta}G_j(z_t)] \over (1-\langle z_t, e_1\rangle)^{1/2}}
&\le \left( {\beta\over 2} + \Re\left[{G_1(z_t)\over1-\langle z_t, e_1\rangle}\right]\left(1 - {t \over 1+t -\eps^2(1-t)^{2\alpha}}\right)
\right){ 1+t -\eps^2(1-t)^{2\alpha}\over \eps(1-t)^\alpha}\cr
&={\beta\over 2}\cdot { 1+t -\eps^2(1-t)^{2\alpha}\over \eps(1-t)^\alpha} + \Re\left[{G_1(z_t)\over1-\langle z_t, e_1\rangle}\right]\left({ 1+t -\eps^2(1-t)^{2\alpha}\over \eps(1-t)^\alpha}  - {t \over \eps(1-t)^\alpha}\right)\cr
&={\beta\over 2}\cdot { 1+t -\eps^2(1-t)^{2\alpha}\over \eps(1-t)^\alpha}+\bigl(-\beta+o\bigl((1-t)^\alpha\bigr)\bigr)
\left({ 1-\eps^2(1-t)^{2\alpha}\over \eps(1-t)^\alpha}\right)\cr
&={\beta\over 2}{\eps^2(1-t)^{\alpha}-(1-t)^{1-\alpha}\over\eps}+o(1)\;.
\cr
%&={\beta\over 2}{\eps^2(1-t)^{\alpha-\gamma+1/2}-(1-t)^{3/2-\alpha-\gamma}\over\eps}+o\left({ 1-\eps^2(1-t)^{2\alpha}\over \eps(1-t)^{\alpha+\gamma-1/2}}\right)
%\cr
}
$$
%Taking $\gamma=1/2-\alpha$ we get
%$$
%{\Re[e^{i\theta}G_j(z_t)] \over (1-\langle z_t, e_1\rangle)^{\gamma}}\le
%{\beta\over 2}{\eps^2(1-t)^{2\alpha}-(1-t)\over\eps}+o(1)
%$$
Letting $t\to 1^-$ we obtain
$$
\limsup_{t\to 1^-}{\Re[e^{i\theta}G_j(z_t)] \over (1-\langle z_t, e_1\rangle)^{1/2}}\le 0
$$
for all $\eps>0$ and $\theta\in\R$. Now letting $\eps\to 0^+$ we find
$$
\limsup_{t\to 1^-}{\Re[e^{i\theta}G_j(te_1)] \over (1-t)^{1/2}}\le 0
$$
for all $\theta\in\R$, and this is possible if and only if
$$
\lim_{t\to 1^-}{G_j(te_1)\over (1-t)^{1/2}}=0\;,
$$
and we are done.
%
%Statements (i), (iii) and (vi), and the $K$-boundedness in (ii), (iv) and (v) were already proved by Bracci and Shoikhet in Theorem \zJWCBS.
%To prove that the restricted $K$-limit in (ii) vanishes, it suffices to argue as in the proof of Theorem \zJWCAR.(ii), but with the following special curve $M$-restricted for any $M>1$: 
%$$
%(0,1)\ni t\mapsto z_t = te_1 + e^{-i\theta} \eps(1-t)^{1/2+\alpha}e_j\in B^n\;,
%$$
%where $j\in\{2,\ldots,n\}$, $0<\eps<1$, $\theta\in\R$, and $\alpha>0$ is the one appearing in the definition of H\"older boundary null point, and we clearly can assume $\alpha<1$. 
%Finally, arguing exactly as in Theorem~\zJWCAR.(iv)--(v), but with $\gamma=1/2$, we obtain (iv) and (v). 
\qedn

We end this paper giving examples of infinitesimal generators having a H\"older boundary null point.

\newex Let $p=e_1$, and $G\colon B^n\to\C^n$ be an infinitesimal generator with 
$\rKlim\limits_{z\to e_1}G(z)=O$. Setting $G_1=\langle G, e_1\rangle$, condition~\eqH\ can be written as
$$
G_1\bigl(\sigma(t)\bigr)=\beta(t-1)+o\bigl((1-t)^{1+\alpha}\bigr)
$$ 
for any special $e_1$-curve $\sigma\colon[0,1)\to B^n$ such that $\langle\sigma(t),e_1\rangle
\equiv t$. In particular, if $G_1$ is of class $C^{1+\alpha'}$ at~$e_1$ for some $\alpha'>\alpha$ then \eqH\ is satisfied, and $e_1$ is a H\"older boundary null point for $G$.  

To give an explicit example, let us recall that if $F\colon B^n\to B^n$ is a holomorphic self-map of~$B^n$ then $G=\id-F$ is an infinitesimal generator (see, e.g., [RS2, Theorem~6.16] and [S, Corollary~3.3.1]). Recalling Theorem~\zJWCR, to get an example of infinitesimal generator having $e_1$ as H\"older boundary null point and satisfying the hypotheses of Theorem~\zJWCAR\ it thus suffices to find $F$ having $K$-limit~$e_1$ at~$e_1$, with $\displaystyle\liminf_{z\to e_1}(1-\|F(z)\|)/(1-\|z\|)<+\infty$ and such that 
$$
F_1\bigl(\sigma(t)\bigr)=t+\beta(1-t)+o\bigl((1-t)^{1+\alpha}\bigr)
$$
for any special $e_1$-curve $\sigma\colon[0,1)\to B^n$ such that $\langle\sigma(t),e_1\rangle\equiv t$. For example, we can just take maps of the form $F(z)=f(z_1) e_1$ with $f$ given by
$$
f(\zeta)=\zeta+\beta(1-\zeta)+c(1-\zeta)^{1+\alpha'}=1-a(1-\zeta)+c(1-\zeta)^{1+\alpha'}\;;
\neweq\eqf
$$
thus we just need to choose $a=1-\beta>0$ and $c>0$ so that $f(\Delta)\subseteq\Delta$. Put $w=1-\zeta$; then $|f(\zeta)|<1$ if and only if $|1-aw+c w^{1+\alpha'}|<1$ if and only if 
$$
a^2|w|^2+c^2|w|^{2(1+\alpha')}+2 c\Re(w^{1+\alpha'})<2a\Re w+2ac|w|^2\Re(w^{\alpha'})\;.
\neweq\eqestu
$$
First of all, write $w=|w|e^{i\theta}$, with $|\theta|<\pi/2$. Then
$$
\Re(w^{\alpha'})=|w|^{\alpha'}\cos(\alpha'\theta)\ge \eps_{\alpha'}|w|^{\alpha'}\;,
$$
where $\eps_{\alpha'}=\cos(\alpha'\pi/2)>0$. Recalling that $|w|<2$, it follows that taking
$c<2^{1-\alpha'}\eps_{\alpha'}a$ we get
$$
c^2|w|^{2(1+\alpha')}<2^{\alpha'}c^2|w|^{2+\alpha'}< 2ac\eps_{\alpha'}|w|^{2+\alpha'}\le
2ac|w|^2\Re(w^{\alpha'})\;.
\neweq\eqstd
$$
Now, if $|\theta|\ge \pi/2(1+\alpha')$ then $\Re(w^{1+\alpha'})\le 0$. Since $|1-w|<1$ implies
$|w|^2<2\Re(w)$, in this case we get
$$
a^2|w|^2+2c\Re(w^{1+\alpha'})<2 a^2\Re w<2 a\Re w
\neweq\eqstt
$$
as soon as $a<1$. 

If instead $|\theta|<\pi/2(1+\alpha')$, we have $|\Im w|<C_{\alpha'}\Re w$, where $C_{\alpha'}=\tan\bigl(\pi/2(1+\alpha')\bigr)$, and thus $|w|<D_{\alpha'}\Re w$, where $D_{\alpha'}=\sqrt{1+C_{\alpha'}^2}$. Hence
$$
a^2|w|^2+2c\Re(w^{1+\alpha'})<\bigl[2 a^2+2cD_{\alpha'}^{1+\alpha'}(\Re w)^{\alpha'}\bigr]\Re w
<2a\Re w
\neweq\eqstq
$$
as soon as $a^2+2^{\alpha'}D_{\alpha'}^{1+\alpha'}c<a$. Since we already requested that $c<2^{1-\alpha'}\eps_{\alpha'}a$, it suffices to have $a<(1+2\eps_{\alpha'}D_{\alpha'}^{1+\alpha'})^{-1}$.

Putting together \eqestu, \eqstd, \eqstt\ and \eqstq\ it follows that if $a<(1+2\eps_{\alpha'}D_{\alpha'}^{1+\alpha'})^{-1}$ and $c<2^{1-\alpha'}\eps_{\alpha'}a$, then the function $f$
given by \eqf\ maps~$\Delta$ into itself, as we wanted.

\setref{EKRS}
\beginsection References

\book A1 M. Abate: Iteration theory of holomorphic maps on taut manifolds!
Me\-di\-ter\-ranean Press, Rende, 1989

\art A2 M. Abate: The Lindel\"of principle and the angular derivative in
strongly convex domains! J. Analyse Math.! 54 1990 189-228

\art A3 M. Abate: Angular derivatives in strongly pseudoconvex domains! Proc. Symp. 
Pure Math.! {52, \rm Part 2} 1991 23-40

\art A4 M. Abate: The infinitesimal generators of semigroups of holomorphic maps! Ann. Mat. Pura Appl.! 161 1992 167-180

\art A5 M. Abate: The Julia-Wolff-Carath\'eodory theorem in polydisks! J.
Analyse Math.! 74 1998 275-306 

\coll A6 M. Abate: Angular derivatives in several complex variables! Real methods in
complex and CR geometry! Eds. D. Zaitsev, G. Zampieri, Lect. Notes in Math.
1848, Springer, Berlin, 2004, pp. 1--47 

\art AT M. Abate,  R. Tauraso: The Lindel\"of principle and angular derivatives
in convex domains of finite type! J. Austr. Math. Soc.! 73 
2002 221-250

\art AMY J. Agler, J.E. McCarthy, N.J. Young: A Carath\'eodory theorem for the bidisk via Hilbert space methods! Math. Ann.! 352 2012 581-624

\art BCD F. Bracci, M.D. Contreras, S. D\'\i az-Madrigal: Pluripotential theory, semigroups and boun\-dary behavior of infinitesimal generators in strongly convex domains! J. Eur. Math. Soc.! 12 2010 23-53

\art BS F. Bracci, D. Shoikhet: Boundary behavior of infinitesimal generators in the unit ball! Trans. Amer. Math. Soc.! 366 no. 2 2014 1119-1140

\book B R.B. Burckel: An introduction to classical complex analysis! Academic Press, New York,1979

\art C C. Carath\'eodory: \"Uber die Winkelderivierten von beschr\"ankten
analytischen Funktionen! Sitzungsber. Preuss. Akad. Wiss. Berlin! {} 1929 39-54

\art {\v C} E.M. \v Cirka: The Lindel\"of and Fatou theorems in $\C^n$! Math.
USSR-Sb.! 21 1973 619-641

\art EKRS M. Elin, D. Khavinson, S. Reich, D. Shoikhet: Linearization models for parabolic dynamical systems via Abel's functional equation! Ann. Acad. Sci. Fen.! 35 2010 1-34

\pre EJ M. Elin, F. Jacobzon: Parabolic type semigroups: asymptotics and order of contact!
Preprint, arXiv:1309.4002! 2013

\art ERS M. Elin, S. Reich, D. Shoikhet: A Julia-Carath\'eodory theorem for hyperbolically monotone mappings in the Hilbert ball! Israel J. Math.! 164 2008 397-411

\art ESY M. Elin, D. Shoikhet, F. Yacobzon: Linearization models for parabolic type semigroups! J. Nonlinear Convex Anal.! 9 2008 205-214

\art H M. Herv\'e: Quelques propri\'et\'es des applications analytiques d'une
boule \`a $m$ dimensions dans elle-m\^eme! J. Math. Pures Appl.! 42 1963 117-147

\art Ju1 G. Julia: M\'emoire sur l'it\'eration des fonctions rationnelles! J. Math. Pures Appl.! 1 1918 47-245

\art Ju2 G. Julia: Extension nouvelle d'un lemme de Schwarz! Acta Math.! 42 1920 349-355

\art Ko A. Kor\'anyi: Harmonic functions on hermitian hyperbolic spaces! Trans.
Amer. Math. Soc.! 135 1969 507-516

\art K-S A. Kor\'anyi, E.M. Stein: Fatou's theorem for generalized
half-planes! Ann. Scuola Norm. Sup. Pisa! 22 1968 107-112

\art L-V E. Landau, G. Valiron: A deduction from Schwarz's lemma! J. London Math. Soc.! 4 1929 162-163

\art Li E. Lindel\"of: Sur un principe g\'en\'erale de l'analyse et ses
applications \`a la theorie de la repr\'esentation conforme! Acta Soc. Sci.
Fennicae! 46 1915 1-35

\art N R. Nevanlinna: Remarques sur le lemme de Schwarz! C.R. Acad. Sci. Paris! 188 1929 1027-1029

\art RS1 S. Reich, D. Shoikhet: Semigroups and generators on convex domains with the hyperbolic metric! Atti Acc. Naz. Lincei Cl. Sc. Fis. Mat. Nat. Rend. Lincei! 8 1997 231-250

\book RS2 S. Reich, D. Shoikhet: Nonlinear semigroups, fixed points, and geometry of domains in Banach spaces! Imperial College Press, London, 2005

\book R W. Rudin: Function theory in the unit ball of $\C^n$! Springer, Berlin, 1980

\book S D. Shoikhet: Semigroups in geometrical function theory! Kluwer Academic Publishers, Dordrecht, 2001

\book St E.M. Stein: The boundary behavior of holomorphic functions of several
complex variables! Princeton University Press, Princeton, 1972

\art Wo J. Wolff: Sur une g\'en\'eralisation d'un th\'eor\`eme de Schwarz! C.R. Acad. Sci. Paris! 183 1926 500-502

\bye